\documentclass[lettersize,journal]{IEEEtran}
\usepackage{amsmath,amsfonts}
\usepackage{algorithmic}
\usepackage{algorithm}
\usepackage{array}
\usepackage[caption=false,font=normalsize,labelfont=sf,textfont=sf]{subfig}
\usepackage{textcomp}
\usepackage{stfloats}
\usepackage{url}
\usepackage{verbatim}
\usepackage{graphicx}
\usepackage{cite}
\usepackage{amssymb}
\usepackage{enumitem}
\usepackage[table]{xcolor} 
\usepackage{colortbl}      
\usepackage{caption}
\usepackage{subcaption}

\newtheorem{proposition}{Proposition}

\usepackage{tabularx}
\usepackage{makecell}
\usepackage{booktabs}
\usepackage{threeparttable}

\begin{document}

\title{Input Convex Neural Network as a Surrogate in Stability-Constrained Optimization for IBR-dominated Power Systems}


\author{Wangkun Xu,~\IEEEmembership{Member,~IEEE},
    Hongyang Jia,~\IEEEmembership{Student Member,~IEEE},
    Yi Wang,~\IEEEmembership{Member,~IEEE}, \\ Ning Zhang,~\IEEEmembership{Senior Member,~IEEE} 
    and Fei Teng,~\IEEEmembership{Senior Member,~IEEE}
    \thanks{
    This is the extended version of ``Input Convex Neural Network as a Surrogate in Stability-Constrained Optimization for IBR-dominated Power Systems''. \\
    This work was partially supported by EPSRC under Grant EP/Y025946/1. \textit{(Corresponding author: Fei Teng)}\\
    W. Xu and F. Teng are with the Department of EEE, Imperial College London, London, UK.
    H. Jia and N. Zhang are with the Department of EE, Tsinghua University, Beijing, China.
    Yi Wang is with the School of Electrical Engineering and Automation, Wuhan University, Wuhan, China.
    Supplementary material including proofs, algorithms, and codes are available at https://github.com/xuwkk/icnn\_constraint\_learning.
    }
}

\markboth{Journal of \LaTeX\ Class Files,~Vol.~14, No.~8, August~2021}%
{Shell \MakeLowercase{\textit{et al.}}: A Sample Article Using IEEEtran.cls for IEEE Journals}


\maketitle

\begin{abstract}
Input convex neural networks (ICNNs) are increasingly used as surrogates for stability indices and embedded as constraints in power-system optimization. This letter clarifies two recurring formulation limitations that can negate ICNN convexity benefits: (i) applying generic Big-$M$ mixed-integer reformulations introduces auxiliary binaries that are unnecessary for enforcing ICNN sublevel constraints; and (ii) reversing the stability inequality transforms a convex sublevel set into a generally nonconvex superlevel set, invalidating global-convergence guarantees of cut-based methods. After clarifying the limitations, we provide (i) an exact LP-based epigraph reformulation for ReLU-ICNNs, (ii) an outer-approximation scheme with global guarantees under the sublevel convention, and (iii) a feasibility-preserving inner-approximation scheme for the superlevel convention, with simulations on IEEE 14- and 118-bus unit commitment instances.
\end{abstract}

\begin{IEEEkeywords}
Power system stability, stability-constrained optimization, constraint learning, input convex neural network.
\end{IEEEkeywords}

\section{Introduction}

The increasing penetration of inverter-based resources (IBRs) has challenged the power system stability and it is urgent to deliver sufficient stability margin at the operational optimization stage. However, due to the complexity of dynamic simulations, real-time stability analysis is difficult to incorporate into optimization for IBR-dominated power systems. Due to their strong representational capability, recent studies have demonstrated the effectiveness of neural-network (NN) surrogate models for approximating power-system stability indices \cite{zhang2022data,xu2025lapso}. Among all classes of NNs, ICNNs \cite{amos2017input} are an appealing candidate. As highlighted in Table~\ref{tab:icnn_pitfalls}, this letter examines two recent studies \cite{wang2025regional,wu2023transient} that do not fully exploit the convexity guarantees offered by ICNNs. 

This letter clarifies two points. (i) From the perspective of a direct (one-shot) formulation \cite{wang2025regional}, an ICNN-based stability constraint can be imposed \emph{without} introducing binary variables. (ii) From the perspective of a cut-based iterative solution \cite{wu2023transient}, guaranteed global convergence can be achieved when the underlying operational problem is (mixed-integer) convex by leveraging convex-analysis arguments and standard decomposition frameworks, including the outer approximation. Finally, we discuss a genuinely nonconvex setting in which the modeling approach of \cite{wu2023transient} remains useful, and we outline how it can be applied appropriately. 

\begin{table}[t]
\caption{Common pitfalls when embedding ICNN constraints}
\label{tab:icnn_pitfalls}
\centering
\setlength{\tabcolsep}{3pt}
\renewcommand{\arraystretch}{1.05}
\begin{tabularx}{\columnwidth}{lXX}
\toprule
& \textbf{\cite{wang2025regional}: Frequency stability} & \textbf{\cite{wu2023transient}: Transient stability} \\
\midrule

\textbf{Solution}
& Direct (one-shot) formulation.
& Iterative formulation. \\\hline

\textbf{Limitation}
& Encoding the ICNN as a generic nonconvex NN leads to auxiliary binary variables that are not required for enforcing the convex ICNN epigraph constraint. 
& Using the reversed inequality ($\geq$) makes the feasible set a \emph{superlevel} set, which is nonconvex; this undermines global convergence guarantees of convex iterative methods. \\
\bottomrule
\end{tabularx}
\end{table}

\section{Preliminaries}

\subsection{Formulation}

This letter focuses on power system operational optimization under a DC power flow model, which can be written as the quadratic program (QP) in \eqref{eq:ed}. Let $z$ denote the decision variables (e.g., generator dispatch) and let $y$ denote input parameters (e.g., renewable generation). To ensure that the resulting operating point remains stable after a disturbance, we embed an additional \emph{algebraic} surrogate constraint $\psi(z;\theta)$, parameterized by $\theta$.
\begin{equation}\label{eq:ed}
    \begin{aligned}
        \textstyle \min_z \; & z^TP(y)z + q^T(y)z \\
        \text{s.t.} \; & A(y)z = b(y),\; G(y)z \leq h(y), \;\psi(z;\theta) \gtreqless \phi_c
    \end{aligned}
\end{equation}
where $P,q,A,G,b,h$ are functions whose structure depends on the grid topology and problem settings. Here $\phi_c$ is a prescribed stability threshold and the sign in $\psi(z;\theta)\,\gtreqless\,\phi_c$ depends on whether the stability metric is defined as ``larger'' or ``smaller-is-more-stable''. It also determines whether the feasible stability set is a sublevel or superlevel set. For convenience, we refer to these two cases as \eqref{eq:ed}-$\leq$ and \eqref{eq:ed}-$\geq$, respectively.

The analysis below is presented for the QP in \eqref{eq:ed}, but the same ideas also apply to more general convex formulations and, under suitable assumptions, to mixed-integer extensions.
Although the proposed framework is general, we focus on a tractable surrogate for small-signal stability \cite{xu2025lapso}, modeled by the generalized short-circuit ratio by \eqref{eq:small_signal_constraint_main_body}. In weak grids, the PCC voltage of an IBR can be strongly influenced by its injected current. This coupling may create a self-synchronizing loop that can degrade the synchronization stability of PLL-based grid-following converters.
\begin{equation}\label{eq:small_signal_constraint_main_body}
    \mathrm{gSCR} = \lambda_{\min}\!\left({Y}_{eq}\right) \geq \mathrm{gSCR}_{\mathrm{lim}}, {Y}_{eq} = \mathrm{diag}\!\left(\frac{{v}_r^{2}}{{p}_r}\right){Y}_{red},
\end{equation}
where ${v}_r$ and ${p}_r$ are the GFL terminal-voltage magnitudes and active-power injections, respectively. The matrix ${Y}_{red}$ is the reduced nodal admittance matrix, as a function of the reactance of online generators. Detailed derivation on gSCR \eqref{eq:small_signal_constraint_main_body} can  be found in Appendix~\ref{app:gscr_detail}.

\subsection{ICNN Architecture}

Due to its convexity on the input under suitable architectural constraints, ICNN becomes a promising candidate surrogate $\psi(z;\theta)$. Consider an ICNN with input $z$ \cite{amos2017input},
\begin{equation}\label{eq:icnn}
    \begin{aligned}
        x_{i+1} & = g_i\bigl(W_i^{x}x_i + W_i^{z}z + b_i\bigr), \; \forall i=0,\cdots,k-1 \\
        \psi(z;\theta) & = x_k,
    \end{aligned}
\end{equation}
where $x_i$ denotes the output of layer $i$ with $x_0 = 0$. $g_{0:k-1}(\cdot)$ are the activation functions. The trainable parameters are denoted by $\theta=\{W_{0:k-1}^{z},W_{1:k-1}^x,b_{0:k-1}\}$ with $W_{0}^{x}=0$. 

\begin{proposition}\label{prop:icnn_convex}
    The ICNN $\psi(z;\theta)$ is convex in $z$ if the following sufficient conditions hold: (i) $W_i^x\geq 0$ elementwise for all $i=1,\cdots,k-1$; (ii) each activation $g_i(\cdot)$ is applied elementwise and is convex for $i=0,\cdots,k-1$; and (iii) $g_i(\cdot)$ is non-decreasing for $i=1,\cdots,k-1$.
\end{proposition}

The proof to this proposation is in Appendix~\ref{app:proof_icnn_convex}.
Compared to the observation in  \cite{amos2017input}, we note that the first activation $g_0$ does not need to be non-decreasing for convexity. The ICNN can also take $y$ as an additional input. In this case, convexity is only required with respect to $z$, yielding a partially input-convex model. When there are multiple stability indices, the ICNN can be designed with vector outputs, where each component is convex in $z$.

\section{ICNN as Stability Constraint}

\subsection{Direct Solution Method}

Encoding conventional NNs with linear layers and piecewise-linear activations (e.g. ReLU) into optimization typically requires auxiliary binary variables, substantially increasing the computational burden of \eqref{eq:ed} \cite{xu2025lapso} (The MILP-based stability-constrained optimization can be found in \eqref{eq:MILP_stability_constraint}). We show that ICNNs can bypass this shortcoming by leveraging convexity.

\begin{proposition}\label{prop:lp_equivalence}
    Under the conditions of Proposition~\ref{prop:icnn_convex}, further assume that ReLU activation is used for $i=0,\cdots,k-2$ and $g_{k-1}(\cdot)$ is the identity map. Then the ICNN \eqref{eq:icnn} is equivalent to the following linear program (LP), in the sense that $\psi(z;\theta)$ equals its optimal value. 
    \begin{equation}\label{eq:lp_equivalence}
    \begin{aligned}
        \psi(z;\theta) & := \textstyle \min_{x_{1:k}} \; x_k \\
        \text{s.t.} \; & x_{i+1} \geq W_i^x x_i + W_i^z z + b_i,\; i=0,\cdots,k-1 \\
        & x_i \geq 0,\; i=1,\cdots,k-1.
    \end{aligned}
    \end{equation}
    Consequently, $\psi(z;\theta)$ is convex, continuous, and piecewise linear in $z$.
\end{proposition}

The proof can be found in Appendix~\ref{app:proof_lp_equivalance}.
Viewing $\psi(z;\theta)$ in Proposition~\ref{prop:lp_equivalence} as an LP value function, the constraint $\psi(z;\theta)\leq \phi_c$ admits the following exact epigraph reformulation.

\begin{proposition}\label{prop:ed_icnn_lp}
    Under the assumption of Proposition~\ref{prop:lp_equivalence}, problem \eqref{eq:ed}-$\leq$ is equivalent to
    \begin{equation}\label{eq:ed_icnn_lp}
    \begin{aligned}
        \textstyle \min_{z,x_{1:k}} \; & z^TP(y)z + q^T(y)z \\
        \text{s.t.} \; & A(y)z = b(y),\; G(y)z \leq h(y), \; x_k \leq \phi_c \\
        & x_{i+1} \geq W_i^x x_i + W_i^z z + b_i,\; i=0,\cdots,k-1 \\
        & x_i \geq 0,\; i=1,\cdots,k-1.
    \end{aligned}
\end{equation}
\end{proposition}

The proof can be found at Appendix~\ref{app:proof_ed_icnn_lp}.
In Proposition~\ref{prop:lp_equivalence}, the ICNN output $\psi(z;\theta)$ equals the optimal value of the LP \eqref{eq:lp_equivalence}. For a fixed $z$, the corresponding optimal primal solution $x_{1:k}$ is not necessarily unique, and the ICNN forward activations provide one such optimal solution. Likewise, in Proposition~\ref{prop:ed_icnn_lp}, the auxiliary variables $x_{1:k}$ need not coincide with the forward-pass activations at the optimizer. Instead, they serve as epigraph variables certifying the existence of a feasible LP solution with $x_k \le \phi_c$. This shows that it is unnecessary to use a generic mixed-integer encoding of ICNNs, as in \cite{wang2025regional} (see Table~\ref{tab:icnn_pitfalls}). Moreover, the same reformulation remains valid when \eqref{eq:ed}-$\leq$ is extended with additional integer variables, since the ICNN contribution is still represented exactly by linear constraints.

\subsection{Cut-based Iterative Solution Method}

\subsubsection{Outer Approximation Method}\label{sec:outer_approximation}

A key benefit of enforcing $\psi(z;\theta)\leq\phi_c$ is that the ICNN output defines a convex sublevel set in $z$, which enables an outer-approximation (OA) algorithm. At iteration $\nu$, the following master problem is solved:
\begin{equation}\label{eq:oa_master}
    \begin{aligned}
        z^{(\nu)} :=\; & \textstyle \arg\min_{z} \; z^TP(y)z + q^T(y)z \\
        \text{s.t.} \; & A(y)z = b(y), \; G(y)z \leq h(y), \\
        & \psi(z^{(s)};\theta) + \lambda^{(s),T} (z - z^{(s)}) \leq \phi_c,\; s=1,\cdots,\nu-1.
    \end{aligned}
\end{equation}
As $\psi(z;\theta)$ is piecewise linear (Proposition~\ref{prop:lp_equivalence}), $\lambda^{(s)} \in \partial \psi(z^{(s)};\theta)$ denotes a subgradient of $\psi$ at $z^{(s)}$, which may be computed via automatic differentiation (e.g., by \texttt{PyTorch}). If the current iterate violates the ICNN constraint, i.e., $\psi(z^{(\nu)};\theta) > \phi_c$, a new cut
\begin{equation}\label{eq:oa_cut}
    \psi(z^{(\nu)};\theta) + \lambda^{(\nu),T} (z - z^{(\nu)}) \leq \phi_c,
\end{equation}
is added.  
Finally, the \emph{global} convergence of OA method is stated as follows.
\begin{proposition}\label{prop:oa_convergence_weak}
    Consider problem \eqref{eq:ed}-$\leq$ under the assumptions of Proposition~\ref{prop:lp_equivalence}. Assume that \eqref{eq:ed}-$\leq$ is feasible and that each master problem \eqref{eq:oa_master} is solved to global optimality. 
    Then: (i) each master problem \eqref{eq:oa_master} is a relaxation of \eqref{eq:ed}-$\leq$, and its optimal value is a lower bound on the optimal value of \eqref{eq:ed}-$\leq$; (ii) the sequence of master optimal values is monotonically non-decreasing and bounded above by the optimal value of \eqref{eq:ed}-$\leq$; (iii) if, at some iteration, the incumbent $z^{(\nu)}$ satisfies $\psi(z^{(\nu)};\theta)\le \phi_c$, then $z^{(\nu)}$ is a global optimum of \eqref{eq:ed}-$\leq$.
\end{proposition}

The proof can be found in Appendix~\ref{app:proof_oa}.
The same statements remain valid for mixed-integer extensions of \eqref{eq:ed}-$\leq$, provided that $\psi(\cdot;\theta)$ remains convex in the relaxed decision vector and each mixed-integer master problem is solved to global optimality.

For multiple ICNN constraints, one may add all violated cuts at each iteration, or only the cut associated with the most violated constraint. In addition, compared with the direct reformulation in Proposition~\ref{prop:ed_icnn_lp}, which introduces auxiliary variables and constraints proportional to the ICNN depth and width, the OA method keeps the master problem smaller and enriches it progressively. Its correctness relies on the convexity of $\psi(\cdot;\theta)$ under the $\leq$ embedding: each OA cut is a valid lower affine approximation of $\psi(\cdot;\theta)$ and hence yields a relaxation of \eqref{eq:ed}-$\leq$. In contrast, under the $\geq$ embedding, the feasible region $\{z: \psi(z;\theta)\ge \phi_c\}$ is generally nonconvex. As a result, the Benders-type decomposition in \cite{wu2023transient} cannot, in general, provide a global-optimality guarantee (see Table~\ref{tab:icnn_pitfalls} and Fig.~\ref{fig:icnn}).

Detailed OA algorithm is provided in Appendix~\ref{app:outer_algorithm}.

\begin{figure}
    \centering
    \includegraphics[width=0.7\linewidth,trim={0 0 0 20 cm},clip]{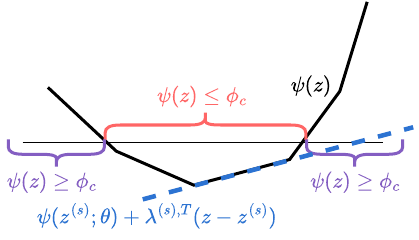}
    \caption{Illustration of ICNN constraint sets $\{z:\psi(z;\theta)\leq \phi_c\}$ and $\{z:\psi(z;\theta)\geq \phi_c\}$ as well as an example of cutting plane.}
    \label{fig:icnn}
\end{figure}

\subsubsection{Inner Approximation Method}\label{sec:inner_approximation}

Although the stability constraint in \cite{wu2023transient} does not define a convex feasible set, the \eqref{eq:ed}-$\geq$ setting can still be handled via an inner-approximation (IA) scheme. We summarize a \emph{local} convergence result below.

Start from a \emph{feasible} point $z^{(0)}$ satisfying $\psi(z^{(0)};\theta)\geq \phi_c$. At iteration $\nu$, compute a subgradient $\lambda^{(\nu-1)}\in\partial \psi(z^{(\nu-1)};\theta)$. Then solve
\begin{equation}\label{eq:inner_approximation}
    \begin{aligned}
        z^{(\nu)} := & \textstyle \arg\min_{z}  z^TP(y)z + q^T(y)z + \textstyle \frac{\rho}{2}\|z - z^{(\nu-1)}\|^2 \\
        \text{s.t.} \; & A(y)z = b(y), \; G(y)z \leq h(y), \\
        & \psi(z^{(\nu-1)};\theta) + \lambda^{(\nu-1),T} (z - z^{(\nu-1)}) \geq \phi_c.
    \end{aligned}
\end{equation}
The proximal term $\frac{\rho}{2}\|z - z^{(\nu-1)}\|^2$ stabilizes the iteration and facilitates convergence.

\begin{proposition}\label{prop:convergence_inner}
    Start from a feasible point $z^{(0)}\in \{z: G(y)z \leq h(y),\; A(y)z = b(y),\; \psi(z;\theta) \geq \phi_c\}$. Under some regularity assumptions, (i) the IA iterates remain feasible, (ii) the objective is non-increasing, and (iii) the method converges to a Karush--Kuhn--Tucker (KKT) stationary point of \eqref{eq:ed}-$\geq$.
\end{proposition}

The proof can be found in Appendix~\ref{app:proof_convergence_inner}. 
A similar conclusion holds for a hybrid mixed-integer implementation in which the integer variables are fixed after a finite number of iterations.

It is essential to start from a feasible point; otherwise, \eqref{eq:inner_approximation} may become infeasible. To obtain such a point, one may solve the following Phase-I feasibility problem iteratively:
\begin{equation}\label{eq:feasibility_problem}
    \begin{aligned}
        (z^{(\nu)}, & s^{(\nu)}) :=\; \arg\min_{z,s} \; s + \textstyle \frac{\rho}{2}\bigl(s - s^{(\nu-1)}\bigr)^2 \\
        \text{s.t.} \; & A(y)z = b(y), \; G(y)z \leq h(y), \; s \geq 0, \\
        & \psi(z^{(\nu-1)};\theta) + \lambda^{(\nu-1),T} (z - z^{(\nu-1)}) + s \geq \phi_c.
    \end{aligned}
\end{equation}
Here, $s$ is a slack variable used to compensate for infeasibility of the current linearized constraint. If the subproblem attains $s^{(\nu)}=0$, then $z^{(\nu)}$ is feasible for the original constraint $\psi(z;\theta)\ge \phi_c$. Thus, \eqref{eq:feasibility_problem} and \eqref{eq:inner_approximation} can be viewed as the Phase-I and Phase-II subproblems of the inner-approximation method, respectively.

Unlike OA, IA cut defines an inner approximation of the true superlevel feasible region $\psi(\cdot;\theta)\geq \phi_c$. This guarantees feasibility preservation, but is conservative. Consequently, even when $\psi(\cdot;\theta)$ is convex, global optimality is generally not guaranteed under the $\geq$ embedding (Fig.~\ref{fig:icnn}). In the presence of multiple ICNN constraints, cuts should be added for all surrogates to preserve feasibility of the original problem.

Detailed IA algorithm is provided in Appendix~\ref{app:inner_algorithm}.

\subsection{Design Principle for Stability Constraints}

In the regression setting, the ICNN surrogate $\psi(\cdot)$ can be trained to predict a scalar stability metric. Some metrics (e.g., small-signal stability margins \eqref{eq:small_signal_constraint_main_body} and transient stability indices) are lower-bounded by a security threshold, whereas others (e.g., frequency nadir violations and RoCoF) are upper-bounded. Accordingly, the stability constraint in \eqref{eq:ed} may be posed with either $\geq$ or $\leq$, depending on the physics definition. However, this physical convention does not restrict the surrogate design. One may redefine the learning target so that the stable region is represented by a convex sublevel set. For example, for a lower-bounded metric with critical threshold $\phi_c$, one may instead predict a signed margin so that stable samples satisfy $\psi(z;\theta) \leq -\phi_c$. In the (binary) classification setting, convexity can be preserved by defining the stable class as the negative label, so that stability corresponds to $\psi(z;\theta) \leq 0$.

Detailed comparison between IA and OA algorithms can be found in Appendix~\ref{app:outer_vs_inner}.

\section{Simulations}\label{sec:simulation}

\begin{table}[t]
    \centering
    \setlength{\tabcolsep}{3.6pt}
    \renewcommand{\arraystretch}{0.9}
    \caption{Simulation Results for 14 Bus system}
    \begin{tabular}{ccccccc}
        \toprule
         &  No.Cons. & No.Vars.(Bin) & Time/$s$ & Cost/k\$ & Vio./\% \\
         \midrule
       Basic & \cellcolor{yellow!15} 2673 & \cellcolor{yellow!15} 1080 (216) & \cellcolor{yellow!15} 0.46 & \cellcolor{yellow!15} 269.2  & \cellcolor{yellow!15} 22.92 \\
       D-MILP-$\leq$ & \cellcolor{blue!15} 9561 & \cellcolor{blue!15} 4584 (1896) & \cellcolor{red!15} 2.13 & \cellcolor{green!15} 277.2 & \cellcolor{green!15} 0.00 \\
       D-LP-$\leq$ & \cellcolor{green!15} 6201 & \cellcolor{green!15} 2904 (216) & \cellcolor{green!15} 0.93 & \cellcolor{green!15} 277.2 & \cellcolor{green!15} 0.00 \\
       I-OA-$\leq$ & \cellcolor{red!15}  2673+$\sum N_\nu$ & \cellcolor{yellow!15} 1080 (216) & \cellcolor{red!15} 1.22 & \cellcolor{green!15} 277.2 & \cellcolor{green!15} 0.00 \\\hline
       D-MILP-$\geq$  & \cellcolor{blue!15} 9561 & \cellcolor{blue!15} 4584 (1896) & \cellcolor{blue!15} 18.62 & \cellcolor{red!15} 279.5 & \cellcolor{green!15} 0.00 \\
       I-IA-$\geq$(I) & \cellcolor{yellow!15} 2697 & \cellcolor{yellow!15} 1080 (216) & \cellcolor{red!15} 1.29 &  \cellcolor{blue!15} 491.9 & \cellcolor{green!15} 0.00 \\
       I-IA-$\geq$(I+II) & \cellcolor{yellow!15} 2697 & \cellcolor{yellow!15} 1080 (216) & \cellcolor{red!15} 2.92 & \cellcolor{red!15} 279.5 & \cellcolor{green!15} 0.00 \\
       \bottomrule
    \end{tabular}
    \label{tab:14_bus}
\end{table}

\begin{table}[t]
    \centering
    \setlength{\tabcolsep}{2.5pt}
    \renewcommand{\arraystretch}{0.9}
    \caption{Simulation Results for 118 Bus system}
    \begin{threeparttable}
    \begin{tabular}{ccccccc}
        \toprule
         &  No.Cons. & No.Vars.(Bin) & Time/$s$ & Cost/k\$ & Vio./\% \\
         \midrule
       Basic & \cellcolor{yellow!15} 25371 & \cellcolor{yellow!15} 10512 (2376) & \cellcolor{yellow!15} 32.63 & \cellcolor{yellow!15} 8872.1 & \cellcolor{yellow!15} 26.67 \\
       D-MILP-$\leq$ & \cellcolor{blue!15} 33435 & \cellcolor{blue!15} 15192 (4056) & \cellcolor{green!15} 47.02 & \cellcolor{green!15} 10371.6 & \cellcolor{green!15} 0.00 \\
       D-LP-$\leq$ & \cellcolor{green!15} 30075 & \cellcolor{green!15} 13512 (2376) & \cellcolor{yellow!15} 30.36 & \cellcolor{green!15} 10370.9 & \cellcolor{green!15} 0.00 \\
       I-OA-$\leq$ & \cellcolor{red!15} 25371+$\sum N_\nu$ & \cellcolor{yellow!15} 10512 (2376) & \cellcolor{red!15} 301.45 & \cellcolor{green!15} 10370.1 & \cellcolor{green!15} 0.00 \\\hline
       D-MILP-$\geq$\tnote{a}  & \cellcolor{blue!15} 33435 & \cellcolor{blue!15} 15192 (4056) & \cellcolor{blue!15} 3611.7 & \cellcolor{black!15} 9507.2 & \cellcolor{red!15} 5.56 \\
       I-IA-$\geq$(I+II)\tnote{a} & \cellcolor{yellow!15} 25395 & \cellcolor{yellow!15} 10512 (2376) & \cellcolor{red!15} 215.53 & \cellcolor{black!15} 10895.8 &  \cellcolor{green!15} 0.00 \\
       I-IA-$\geq$(I) & \cellcolor{yellow!15} 25395 & \cellcolor{yellow!15} 10512 (2376) & \cellcolor{green!15} 61.78 & \cellcolor{blue!15} 77839.6 & \cellcolor{red!15} 0.42 \\
       I-IA-$\geq$(I+II) & \cellcolor{yellow!15} 25395 & \cellcolor{yellow!15} 10512 (2376) & \cellcolor{red!15} 269.84 & \cellcolor{red!15} 12632.1 &  \cellcolor{red!15} 2.50 \\
       \bottomrule
    \end{tabular}
    \begin{tablenotes}
        \footnotesize
        \item[a] D-MILP-$\geq$ finds feasible solutions for only 70\% of samples; therefore, D-MILP-$\geq$ and I-IA-$\geq$ here are reported on these feasible cases only.
    \end{tablenotes}
    \end{threeparttable}
    \label{tab:118_bus}
\end{table}

We use realistic load and solar generation profiles over one year and preprocess/augment the training data for small-signal stability assessment into 78{,}840 samples using the procedure in \cite{xu2025lapso}. A full unit commitment with a 24-hour horizon is solved for the IEEE 14-bus and 118-bus systems, where small-signal stability is enforced through an ICNN surrogate with two hidden layers of sizes 50 and 20. More details on the ICNN structure and training loss can be found in Appendix~\ref{app:ICNN_design}. The exact (stability-constrained) UC formulation can be found in Appendix~\ref{app:uc_detail}. Detailed simulation settings, including dataset generation and training hyperparameters can be found in Appendix~\ref{app:simulation_setting}.

Five algorithms are compared: (i) \underline{D-MILP-$\geq$}, direct solution with MILP reformulation under the ``$\geq$'' embedding; (ii) \underline{D-MILP-$\leq$}, direct solution with a MILP reformulation under the ``$\leq$'' embedding (adapted from \cite{wang2025regional}); (iii) \underline{D-LP-$\leq$}, direct solution with the LP reformulation in Proposition~\ref{prop:ed_icnn_lp}; (iv) \underline{I-OA-$\leq$}, iterative solution via outer approximation in Proposition~\ref{prop:oa_convergence_weak}; and (v) \underline{I-IA-$\geq$}, iterative solution via inner approximation in Proposition~\ref{prop:convergence_inner} (adapted from \cite{wu2023transient}). We also report the performance of the basic unit commitment problem without stability constraints. ``I+II'' denotes that both Phase-I \eqref{eq:feasibility_problem} and Phase-II \eqref{eq:inner_approximation} are included while `'I'' represents the Phase-I performance alone.
Tables~\ref{tab:14_bus} and \ref{tab:118_bus} summarize the results, averaged over 10 evenly sampled days across the year; similar performances are highlighted using the same color. Here, $N_\nu$ denotes the number of violated constraints at iteration $\nu$, and \emph{Vio./\%} denotes the percentage of time periods whose operating points are identified as unstable by the \emph{model-based} small-signal stability assessment \eqref{eq:small_signal_constraint_main_body}.

The detailed formulation on D-MILP-$\geq$ and D-MILP-$\leq$ can be found in Appendix~\ref{app:ICNN_mixed_integer}. Detailed OA and IA algorithms for solving ICNN-based stability-constrained UC can be found in Appendix~\ref{app:outer_algorithm} and \ref{app:inner_algorithm}, respectively. Detailed number of constraints and number of extra parameters can be found in Appendix~\ref{app:optimization_structure}.

For both systems, D-LP-$\leq$ achieves the lowest computational time and performs comparably to the basic setting. By contrast, the computational time of D-MILP-$\leq$ is roughly twice as large. In this case study, I-OA-$\leq$ does not provide a clear computational advantage, mainly because the overall solution burden is dominated by the integer variables in the underlying UC, while the number of ICNN-induced constraints is relatively limited. Under the ``$\geq$'' embedding, I-IA-$\geq$ can serve as a practical alternative to the direct MILP formulation. In the 14-bus case, it attains the same objective value as D-MILP-$\geq$, suggesting that the solution of the basic problem provides an effective warm start for the Phase-I problem \eqref{eq:feasibility_problem}. Moreover, all methods substantially improve system stability.

The 118-bus results under the ``$\geq$'' embedding show a different trend. In particular, D-MILP-$\geq$ may fail to find a feasible solution within the one-hour time limit. Moreover, the ICNN surrogate under the ``$\geq$'' embedding achieves only about 90\% prediction accuracy, indicating that the small-signal stability relation in \eqref{eq:small_signal_constraint_main_body} is not well matched to a superlevel-set approximation. By contrast, the stronger performance of the ``$\leq$'' embedding suggests that the same relation \eqref{eq:small_signal_constraint_main_body} is better captured through a convex sublevel-set surrogate. As a result, the ``$\geq$'' embedding requires a more conservative threshold $\phi_c$ to maintain system security, which increases operating cost. In addition, because the I-IA-$\geq$ method only provides local stationarity guarantees, it terminates at a suboptimal local solution, further increasing the operating cost relative to D-MILP-$\geq$ in the 118-bus case.

\section{Conclusion}

This letter formalizes ICNN-based stability constraints for IBR-dominated power system operations and clarifies common formulation limitations that can increase computational burden and undermine convergence. We demonstrate that when stability is encoded through a convex sublevel-set constraint, the ICNN can be exploited exactly through an LP reformulation, making D-LP-$\leq$ more attractive than generic MILP-based encodings. In contrast, the ``$\geq$'' embedding generally induces a nonconvex feasible region, for which MILP or conservative inner-approximation type schemes are available; the inner-approximation preserve feasibility, but does not provide global-optimality guarantees.


\appendix

\subsection{Open Source Development}\label{app:open_source}

All experiment in the paper is open-source at \url{https://github.com/xuwkk/icnn_constraint_learning}. The data generation and optimization formulation is implemented by our open source package \texttt{GridForge}: \url{https://github.com/xuwkk/gridforge}.

\subsection{gSCR-based Small Signal Stability}\label{app:gscr_detail}

\subsubsection{Stability Index}

In weak grids (i.e., with large network impedance), the PCC voltage of an inverter-based resource (IBR) can be strongly influenced by its injected current. This coupling may create a self-synchronizing loop between the PLL dynamics and the network response, effectively introducing a positive-feedback path that can degrade the synchronization stability of PLL-based grid-following (GFL) converters. Following \cite{dong2019small}, \cite{liu2023generalized}, and \cite{xu2025incorporation}, we quantify the small-signal synchronization stability of GFL IBRs using the generalized short-circuit ratio (gSCR). The underlying premise is that synchronization stability is primarily governed by the interaction between converter dynamics and the network dynamics, the latter being captured compactly by the nodal admittance matrix.

Specifically, gSCR is defined through an equivalent network matrix ${Y}_{eq}$:
\begin{equation}\label{eq:small_signal_constraint_1}
\begin{aligned}
    \mathrm{gSCR} &= \lambda_{\min}\!\left({Y}_{eq}\right),\\
    {Y}_{eq} &= \mathrm{diag}\!\left(\frac{{v}_r^{2}}{{p}_r}\right){Y}_{red},
\end{aligned}
\end{equation}
where $\mathrm{diag}\!\left({{v}_r^{2}}/{{p}_r}\right)$ is a diagonal scaling constructed from the GFL terminal-voltage magnitudes ${v}_r$ and active-power injections ${p}_r$ (with ${v}_r^{2}/{p}_r$ interpreted element-wise). The matrix ${Y}_{red}$ is the reduced nodal admittance matrix obtained after eliminating passive buses and infinite buses. The matrix ${Y}_{eq}$ is diagonalizable, and its smallest eigenvalue $\lambda_{\min}({Y}_{eq})\in\mathbb{R}_{+}$ serves as an index of effective network strength seen by the GFL units. Hence, the small-signal synchronization requirement can be enforced using the threshold form \cite{liu2023generalized}
\begin{equation}\label{eq:gSCR_cstrt}
    \mathrm{gSCR} \ge \mathrm{gSCR}_{\mathrm{lim}},
\end{equation}
where $\mathrm{gSCR}_{\mathrm{lim}}$ denotes the critical (minimum) gSCR required for small-signal stability of the GFL units. Moreover, under the standard operating assumption that bus voltages remain close to $1~\mathrm{p.u.}$ for normal conditions and sufficiently small disturbances, $\mathrm{gSCR}_{\mathrm{lim}}$ can be treated as (approximately) operation-independent and thus determined offline, with or without detailed converter-control parameter availability \cite{liu2023generalized}. 

Considering the ICNN surrogate constraints \eqref{eq:ed}, \eqref{eq:gSCR_cstrt} represents a $\geq$ setting ($\psi(z;\theta) \geq \phi_c$). However, it is still possible to learning a \eqref{eq:ed}-$\leq$ surrogate as discussed in Section~\ref{sec:simulation}.

\subsubsection{Constraint Formulation}

As discussed above, $\mathrm{gSCR}_{\mathrm{lim}}$ in \eqref{eq:gSCR_cstrt} is fixed at the scheduling stage (offline tuning/assessment), whereas $\mathrm{gSCR}$ varies with the operating point, in particular the renewable injections and generator commitment decisions, i.e.,
$\mathrm{gSCR}=\psi(z;\theta) = \mathrm{gSCR}({u}_g, P_{rn})$ where $u_g$ and $P_{rn}$ are the generator commitment status and renewable generation, respectively.
We now make this dependence explicit.

Let ${Y}_0$ denote the \emph{nodal admittance matrix} of the transmission network with purely reactive lines \cite{liu2023generalized}. To account for SGs, we consider the augmented admittance
\begin{equation}\label{eq:Y1}
    {Y}={Y}_0+{Y}_g .
\end{equation}
Under normal operation and small disturbances, each SG can be modeled as a voltage source behind a reactance; consequently, the SG contribution ${Y}_g$ depends on the unit commitment and can be written as
\begin{equation}\label{eq:Y2}
    {Y}_{g,ij}=
    \begin{cases}
    \dfrac{1}{x_{g,i}}\,u_{g,i}, & \text{if } i=j\in\mathcal{N}_g,\\[6pt]
    0, & \text{otherwise},
    \end{cases}
\end{equation}
where $x_{g,i}$ is the (internal) reactance of SG $i$, $\mathcal{N}_g$ is the set of SG buses, and $u_{g,i}\in\{0,1\}$ indicates whether unit $i$ is committed.

Next, permute ${Y}$ into the block form
\begin{equation}\label{eq:Y0}
    {Y}=
    \begin{bmatrix}
        {Y}_{\mathcal{N}_r\mathcal{N}_r} & {Y}_{\mathcal{N}_r\delta}\\
        {Y}_{\delta\mathcal{N}_r} & {Y}_{\delta\delta}
    \end{bmatrix},
\end{equation}
where $\mathcal{N}_r\subseteq\mathcal{N}$ is the set of GFL IBR buses and $\delta=\mathcal{N}\setminus\mathcal{N}_r$ collects all remaining buses. The reduced admittance matrix in \eqref{eq:small_signal_constraint_1} is then obtained via Kron reduction:
\begin{equation}\label{Y_red}
    {Y}_{red}
    ={Y}_{\mathcal{N}_r\mathcal{N}_r}
    -{Y}_{\mathcal{N}_r\delta}{Y}_{\delta\delta}^{-1}{Y}_{\delta\mathcal{N}_r}.
\end{equation}
Substituting \eqref{Y_red} into \eqref{eq:small_signal_constraint_1} yields an explicit mapping from the scheduling decisions $({u}_g, P_{rn})$ to ${Y}_{eq}$ and hence to $\mathrm{gSCR}$, enabling enforcement of \eqref{eq:gSCR_cstrt} within the system scheduling model.

\subsection{Proofs on the Propositions}\label{app:proofs}

\subsubsection{Proof of Proposition \ref{prop:icnn_convex}}\label{app:proof_icnn_convex}

As stated in \cite{amos2017input}, the proof directly follows from the composition of convex functions. We prove by the property of element-wise of convex function with induction.

For each element $r$ of the first layer, the ICNN computes,
\begin{equation*}
    [x_1(z)]_r = g_0(\langle [W_0^{z}]_r, z \rangle + [b_0]_r)
\end{equation*}
where $[\cdot]_r$ represents the $r$-th row of a matrix or $r$-th element of a vector.
As the composition of \emph{convex function} $g_0(\cdot)$ with affine mapping is convex \cite{boyd2004convex}, each element of $x_1(z)$ is convex function of $z$. Then starting from the second layer, each output $r$ computes,
\begin{equation*}
    \begin{aligned}
         [s_{i+1}(z)]_r = & \langle [W_i^x]_r, x_i(z) \rangle + \langle [W_i^z]_r, z \rangle + [b_i]_r  \\
          = & \textstyle \sum_j [W_i^x]_{rj} [x_i(z)]_j + \langle [W_i^z]_r, z \rangle + [b_i]_r
    \end{aligned}
\end{equation*}
If $W_i^x\geq 0$ elementwise, then $\sum_j [W_i^x]_{rj} [x_i(z)]_j$ is a nonnegative weighted sum of convex functions and hence convex; adding an affine term preserves convexity, so $[s_{i+1}(z)]_r$ is convex. Since $g_i$ is convex and non-decreasing (for $i\geq 1$), the composition $[x_{i+1}(z)]_r = g_i([x_{i+1}(z)]_r)$ is convex \cite{boyd2004convex}. Therefore each component of $x_{i+1}$ is convex on $z$, completing the induction and proving convexity of $\psi(z;\theta)$. 

\subsubsection{Proof of Proposition \ref{prop:lp_equivalence}}\label{app:proof_lp_equivalance}

Let $\hat{x}_{i},\; i=1,\cdots,k$ the output of ICNN \eqref{eq:icnn}. Let $\tilde{x}_i,\;i=1,\cdots,k$ and $x_i^\star,\;i=1,\cdots,k$ be the feasible point and the optimal solution of LP \eqref{eq:lp_equivalence}, respectively. 

Forward passing $z$ by ICNN \eqref{eq:icnn} with $\hat{x}_{i+1} = \operatorname{ReLU}(W_i^x\hat{x}_i + W_i^zz + b_i)$ gives that 
\begin{equation*}
    \begin{aligned}
        \hat{x}_{i+1} & \geq W_i^x\hat{x}_i + W_i^zz + b_i,~~ i=0,\cdots,k-2 \\ \hat{x}_{i+1} & \geq 0,~~ i=0,\cdots,k-2 \\
        \hat{x}_k & = W_{k-1}^x\hat{x}_k + W_{k-1}^z z + b_{k-1}
    \end{aligned}
\end{equation*}
which is feasible for LP \eqref{eq:lp_equivalence}. Therefore, $x_k^\star \leq \hat{x}_k$.

For LP-feasible $\tilde{x}_1$, it satisfies $\tilde{x}_1\geq W_0^z z + b_0$ and $\tilde{x}_1 \geq 0$. Therefore, 
\begin{equation*}
    \tilde{x}_1 \geq \max\{0,W_0^zz + b_0\} = \operatorname{ReLU}(W_0^zz + b_0) = \hat{x}_1
\end{equation*}
For $i=1,\cdots,k-2$, the feasibility also gives that
\begin{equation*}
        \tilde{x}_{i+1} \geq \operatorname{ReLU}(W_i^x\tilde{x}_i + W_i^zz + b_i)
\end{equation*}
Then by induction, given $\tilde{x}_{1} \geq \hat{x}_1$ and $W_i^{x} \geq 0$, for $i=1,\cdots,k-2$,
\begin{equation*}
    \begin{aligned}
        & W_i^x\tilde{x}_i + W_i^zz + b_i \geq W_i^x\hat{x}_i + W_i^zz + b_i \\
        \Rightarrow \; & \tilde{x}_{i+1} \geq \operatorname{ReLU}(W_i^x\tilde{x}_i + W_i^zz + b_i) \\
        & \hspace{2cm} \geq \operatorname{ReLU}(W_i^x\hat{x}_i + W_i^zz + b_i) = \hat{x}_{i+1}
    \end{aligned}
\end{equation*}
Similarly, for $i=k-1$, the feasibility of \eqref{eq:lp_equivalence} gives,
\begin{equation*}
    \begin{aligned}
        \tilde{x}_{k} & \geq W_{k-1}^x\tilde{x}_{k-1} + W_{k-1}^zz + b_{k-1} \\
        & \hspace{2cm} \geq W_{k-1}^x\hat{x}_{k-1} + W_{k-1}^zz + b_{k-1} = \hat{x}_{k}
    \end{aligned}
\end{equation*}
Therefore, it is proved that $\tilde{x}_i \geq \hat{x}_i,\; i=1,\cdots,k$. Note that this is valid for any feasible $\tilde{x}_i$ including the optimum $x^\star_i$. I.e., $x^\star_k \geq \hat{x}_k$. Together with $x^\star_k \leq \hat{x}_k$, we conclude $x^\star_k = \hat{x}_k$. 

The proof on the continuous and piece-wise linear property directly follows from \cite{pistikopoulos2007multi}, as the ICNN is represented as a parametric LP over $x$.

\subsubsection{Proof of Proposition \ref{prop:ed_icnn_lp}}\label{app:proof_ed_icnn_lp}

Define feasible set of \eqref{eq:lp_equivalence} as $\mathcal{X}(z)=\{x_{1:k}\mid x_{i+1} \geq W_i^x x_i + W_i^z z + b_i,\; i=0,\cdots,k-1,  x_i \geq 0,\; i=1,\cdots,k-1\}$. Then for any fixed $z$, it is easy to show that
\begin{equation}\label{eq:sup_1}
    \psi(z;\theta) \leq \phi_c \Leftrightarrow \exists x_{1:k} \in\mathcal{X}(z) \text{ such that } x_k\leq \phi_c
\end{equation}

Consider a feasible $z$ of \eqref{eq:ed}-$\leq$. Then by \eqref{eq:sup_1}, there exists $x_{1:k}$ such that $(z,x_{1:k})$ is feasible for \eqref{eq:ed_icnn_lp}. Similarly, for a feasible $(z,x_{1:k})$ of \eqref{eq:ed_icnn_lp}, \eqref{eq:sup_1} shows that $\psi(z;\theta) \leq \phi_c$. Therefore, $z$ is feasible for \eqref{eq:ed}-$\leq$. Therefore, the two problems have the same feasible $z$-set, and since both objectives are identical and depend only on $z$, they have the same optimal objective value and the same set of optimal 
$z$-solutions (with $x$ acting as an auxiliary variables in \eqref{eq:ed_icnn_lp}). This proves equivalence.

\subsubsection{Proof of Proposition~\ref{prop:oa_convergence_weak}}
\label{app:proof_oa}
For any previously generated point $z^{(s)}$ and any subgradient $\lambda^{(s)}\in\partial_z\psi(z^{(s)};\theta)$, convexity of $\psi(\cdot;\theta)$ implies
\begin{equation}
    \psi(z;\theta)\ge \psi(z^{(s)};\theta)+\lambda^{(s),T}(z-z^{(s)}),\quad \forall z.
\end{equation}
Hence, any $z$ feasible for \eqref{eq:ed}-$\leq$ ($\psi(z;\theta) \leq \phi_c$) also satisfies every OA cut in \eqref{eq:oa_master}. Therefore, each master problem is a relaxation of \eqref{eq:ed}-$\leq$, and its optimal value is a lower bound on the optimal value of \eqref{eq:ed}-$\leq$. This proves (i).

Since each iteration adds one additional valid constraint to the master problem, the feasible region of the master problem can only shrink. Because each master problem is solved to global optimality, its optimal value is therefore monotonically nondecreasing. By part (i), these optimal values are always lower bounds on the optimal value of \eqref{eq:ed}-$\leq$, and hence they are bounded above by that optimal value. This proves (ii).

Finally, suppose that at some iteration $\nu$, the incumbent $z^{(\nu)}$ satisfies $\psi(z^{(\nu)};\theta)\le \phi_c$. Then $z^{(\nu)}$ is feasible for \eqref{eq:ed}-$\leq$. Since it is also optimal for the current master problem, whose optimal value is a lower bound on the optimal value of \eqref{eq:ed}-$\leq$, the objective value attained by $z^{(\nu)}$ must coincide with the global optimum of \eqref{eq:ed}-$\leq$. Hence $z^{(\nu)}$ is a global optimum. This proves (iii).

The mixed-integer case follows by the same argument, since the proof only uses validity of the subgradient cuts and global solution of each master problem.

\subsubsection{Proof of Proposition~\ref{prop:convergence_inner}}
\label{app:proof_convergence_inner}

Consider the following abstract form of the original problem:
\begin{equation}\label{eq:inner_1}
    \min_{z\in\mathcal{Z}} f(z) \quad \text{s.t.} \; \psi(z;\theta) \geq \phi_c
\end{equation}
where $\mathcal{Z}$ is nonempty, closed, and convex (polyhedral), $f$ is convex and continuously differentiable, $\psi(\cdot)$ is convex and piece-wise linear, and $\rho\geq0$. In addition, assume that \eqref{eq:inner_1} is bounded $f_{\inf} := \inf_{z\in\mathcal{Z},\psi(z;\theta)\geq\phi_c} > -\infty$.

Given a feasible initial point $z^{(0)} \in \mathcal{Z}$ satisfying $\psi(z^{(0)};\theta) \geq \phi_c$, select a subgradient $\lambda^{(s)} \in \partial \psi(z^{(s)};\theta)$ and define the supporting hyperplane at iteration $s$ by
\begin{equation*}
    \ell^{(s)}(z) := \psi(z^{(s)};\theta) + \lambda^{(s),T} (z - z^{(s)})
\end{equation*}
At iteration $s$, the convex subproblem of \eqref{eq:inner_approximation} is compactly written as,
\begin{equation}\label{eq:inner_2}
    z^{(s+1)} = \arg\min_{z\in\mathcal{Z}} (f(z) + \textstyle \frac{\rho}{2} \|z - z^{(s)}\|^2) \quad \text{s.t.} \; \ell^{(s)}(z) \geq \phi_c
\end{equation}

\emph{Feasibility.} Because $\psi(\cdot;\theta)$ is convex, the supporting hyperplane $\ell^{(s)}$ is a global underestimator, i.e., $\psi(z;\theta) \geq \ell^{(s)}(z),\; \forall z$. Therefore, the subproblem constraint $\ell^{(s)}(z) \geq \phi_c$ implies $\psi(z;\theta) \geq \phi_c$. Hence any solution $z^{(s+1)}$ of \eqref{eq:inner_2} satisfies $\psi(z^{(s+1)};\theta) \geq \phi_c$ and $z^{(s+1)} \in \mathcal{Z}$. Since $z^{(0)}$ is feasible, induction shows that $z^{(s)}$ is feasible for every $s\geq0$.

\emph{Monotone non-increasing and convergence.} Note that $z^{(s)}$ is naturally feasible for \eqref{eq:inner_2} as $\ell^{(s)}(z^{(s)}) = \psi(z^{(s)};\theta) \geq \phi_c$. Optimality of $z^{(s+1)}$ for \eqref{eq:inner_2} gives,
\begin{equation}\label{eq:inner_3}
    \begin{aligned}
        & f(z^{(s+1)}) + \textstyle \frac{\rho}{2} \|z^{(s+1)} - z^{(s)}\|^2 \leq f(z^{(s)}) + \textstyle \frac{\rho}{2} \|z^{(s)} - z^{(s)}\|^2 \\
        & \hspace{6cm} = f(z^{(s)})
    \end{aligned}
\end{equation}
which means $f(z^{(s+1)}) \leq f(z^{(s)})$ and hence $\{f(z^{(s)})\}$ is non-increasing.

Summing \eqref{eq:inner_3} over $s=0,\cdots,\nu$ gives,
\begin{equation*}
    \begin{aligned}
        \textstyle \frac{\rho}{2} \sum_{s=0}^{\nu} \|z^{(s+1)} - z^{(s)}\|^2 &\leq f( z^{(0)}) - f( z^{(\nu+1)} ) \\ & \leq f(z^{(0)}) - \textstyle f_{\inf} < \infty
    \end{aligned}
\end{equation*}
The second and third inequality follows from the feasibility of $z^{(\nu+1)}$ for the original problem and the definition of $f_{\inf}$. Letting $\nu\rightarrow\infty$ yields $\sum_{s=0}^{\infty}\|z^{(s+1)}-s^{(s)}\|^2<\infty$. As each term is nonnegative, the terms must go to zero, i.e., $\|z^{(s+1)} - z^{(s)}\| \rightarrow 0$.

\emph{KKT stationary.} Since $\mathcal{Z}$ is polyhedral and the constraint $\ell^{(s)}(z)\geq\phi_c$ is affine, the KKT conditions hold at the solution $z^{(s+1)}$ of \eqref{eq:inner_2}. In particular, there exists a multiplier $\varsigma^{(s+1)} \geq 0$ such that
\begin{equation}\label{eq:inner_4}
    \begin{aligned}
        & \nabla f(z^{(s+1)}) + \rho\cdot (z^{(s+1)} - z^{(s)} ) - \varsigma^{(s+1)} \cdot \lambda^{(s)} \\
        & \hspace{5cm} + N_{\mathcal{Z}}(z^{(s+1)}) \ni 0
    \end{aligned}
\end{equation}
where $N_{\mathcal{Z}}(\cdot)$ is the normal cone of $\mathcal{Z}$. Moreover, the primal feasibility and complementary slackness conditions are
\begin{equation}\label{eq:inner_5}
    \ell^{(s)}(z^{(s+1)})\geq \phi_c, \quad \varsigma^{(s+1)} \cdot (\ell^{(s)}(z^{(s+1)}) - \phi_c) = 0
\end{equation}

We next show that every accumulation point of the sequence $\{z^{(s)}\}$ satisfies the generalized KKT conditions of \eqref{eq:inner_1}. Let $z^\star$ be an arbitrary accumulation point. As $\{z^{(s)}\}$ is bounded, there exists a subsequence $\{s_j\}$ such that $z^{(s_j)}\rightarrow z^\star$. Since $\|z^{(s+1)} - z^{(s)}\|\rightarrow 0$, it also follows that $z^{(s_j+1)}\rightarrow z^\star$. In addition, since $\psi(\cdot;\theta)$ is finite, convex and piecewise linear, its subgradients are bounded. Hence, $\{\lambda^{(s_j)}\}$ has a convergent subsequence. Suppose that $\lambda^{(s_j)}\rightarrow \lambda^\star$ and $\varsigma^{(s_j+1)}\rightarrow\sigma^\star\geq0$. 

For every $z$,
\begin{equation*}
    \psi(z;\theta) \geq \psi(z^{(s_j)};\theta) + \lambda^{(s_j),T}(z-z^{(s_j)})
\end{equation*}
Taking $j\rightarrow\infty$ gives
\begin{equation*}
    \psi(z;\theta) \geq \psi(z^\star;\theta) + \lambda^{\star,T}(z-z^{\star}), ~~ \forall z
\end{equation*}
and therefore $\lambda^\star\in\partial\psi(z^\star;\theta)$; i.e., the converged $\lambda^{(s_j)}$ is the subgradient of $\psi(z^\star;\theta)$. 
Moreover,
\begin{equation}\label{eq:inner_6}
    \ell^{(s_j)}(z^{(s_j+1)}) = \psi(z^{(s_j)};\theta) + \lambda^{(s_j),T}(z^{(s_j+1)} - z^{(s_j)}) \rightarrow \psi(z^\star;\theta)
\end{equation}
where the limit follows from the boundednesss of $\{\lambda^{(s_j)}\}$ and $z^{(s_j+1)}-z^{(s_j)}\rightarrow 0$.

Passing the limit in \eqref{eq:inner_4} gives
\begin{equation*}
    0\in\nabla f(z^\star) - \varsigma^\star\cdot\lambda^\star + N_\mathcal{Z}(z^\star)
\end{equation*}
Similarly, passing the limit \eqref{eq:inner_6} to \eqref{eq:inner_5} yields
\begin{equation*}
    \psi(z^\star;\theta) \geq \phi_c, \quad \varsigma^\star \cdot (\psi(z^\star;\theta) - \phi_c) = 0 
\end{equation*}

Since $z^\star$ was chosen arbitrarily, every accumulation point of the IA sequence is a generalized KKT stationary of \eqref{eq:inner_1}. If the sequence has a unique accumulation point, then the entire sequence converges to that stationary point.

\subsection{Algorithm Design and Properties}\label{app:algorithm}

\subsubsection{ICNN Design and Training}\label{app:ICNN_design}

The ICNN is designed by \eqref{eq:icnn} where the output activation is identity. A `SoftPlus' function is added to keep the $W^x$ nonnegative, which is defined as
\begin{equation*}
    \operatorname{softplus}([W^x]_{ij}) = \log(1+e^{[W^x]_{ij}})
\end{equation*}
The binary cross entropy loss is considered,
\begin{equation*}
    \begin{aligned}
        & \mathcal{L}(\theta,\mathcal{D}) = \textstyle -\frac{1}{N} \sum_{i=1}^N y(i)\cdot\log\sigma(u_g(i),P_{rn}(i),\theta) \\
        & \hspace{2cm} + (1-y(i))\cdot\log(1-\sigma(\psi(u_g(i),P_{rn}(i),\theta))) 
    \end{aligned}
\end{equation*}
where $\sigma(\cdot)$ is the sigmoid function. We implement this as \texttt{BCEWithLogitsLoss} using logits $\psi(z;\theta)$ to avoid numerical underflow.

Moreover, as the UC horizon is 24 hours, the ICNN is trained to assess the small-signal stability of each instance. I.e, 24 ICNN surrogate constraints $\psi(z_t,\theta) \gtreqless 0$ are added in the small-signal stability constrained UC \eqref{eq:ed}.
\begin{equation}\label{eq:ed_over_time}
    \begin{aligned}
         \min_{z_{1:T}} \; & z^TP(y)z + q^T(y)z \\
        \text{s.t.} \; & A(y)z = b(y),\; G(y)z \leq h(y) \\
        & \psi(z_t;\theta) \gtreqless \phi_c,\; t=1,\cdots,T
    \end{aligned}
\end{equation}

\subsubsection{UC Formulation with ICNN Constraints}\label{app:uc_detail}

A standard unit commitment (UC) problem is modified from \cite{conejo2018power}. The parameter and decision variable definitions are reported in Table~\ref{tab:param_uc_paper} and Table~\ref{tab:var_uc_paper}.

The objective \eqref{eq:objective} includes the generator dispatch, online, start-up and shut-down cost, as well as the load shedding and renewable generation curtailment cost. \eqref{eq:commitment_logic_1}-\eqref{eq:commitment_logic_2} enforce the generator commitment logic. Generator ramp-up, ramp-down, and generation limit are constrained by \eqref{eq:ramp_up}-\eqref{eq:gen_limit}. The nodal power injection is considered in \eqref{eq:injection} with power balance and power flow limits in \eqref{eq:power_balance} and \eqref{eq:power_flow}. The non-negativity of slack variables are enforced in \eqref{eq:slack_1}-\eqref{eq:slack_2}.
The generator minimum online and offline time is considered in \eqref{eq:minimin_on}-\eqref{eq:minimum_off}. The ICNN surrogate \eqref{eq:icnn_constratint} of \eqref{eq:small_signal_constraint_1} takes the generator commitment status $u_g(t)$ and the net renewable generation $P_{rn} = P_r(t) - P_{rc}(t)$ as input. Moreover, it is considered that $u_g(0) = u_{g,init}$ and $P_g(0)=P_{g,init}$ are the known initial conditions of this operational period. The set of decision variable is denoted as $\Xi = \{P_g(t),u_g(t),y_g(t),z_g(t), P_{ls}(t), P_{rc}(t)\}_{t=1}^T$.

\begin{subequations}
    \begin{align}
        \min_{\Xi} \; & \textstyle \sum_{t=1}^T c_g^TP_g(t) + c_u^Tu_g(t) + c_y^Ty_g(t) + c_z^Tz_g(t) \nonumber \\
        & \hspace{1.5cm} + c_{ls}^TP_{ls}(t) + c_{rc}P_{rc}(t)  \label{eq:objective}  \\
        \text{s.t.} \; & y_g(t) - z_g(t) = u_g(t) - u_g(t-1), \; t=1,\cdots,T \label{eq:commitment_logic_1} \\ 
        & y_g(t) + z_g(t) \leq 1, \; t=1,\cdots,T  \label{eq:commitment_logic_2} \\
        & P_g(t) - P_g(t-1) \leq R_{up} \circ u_g(t-1) \nonumber \\ 
        & \hspace{1.5cm} + R_{start} \circ y_g(t), \; t=1,\cdots,T \label{eq:ramp_up} \\
        & P_g(t-1) - P_g(t) \leq R_{down} \circ u_g(t) \nonumber\\
        & \hspace{1.5cm} + R_{shut} \circ z_g(t), \; t=1,\cdots,T \label{eq:ramp_down} \\
        & P_{min} \circ u_g(t) \leq P_g(t) \leq P_{max} \circ u_g(t), \; t=1,\cdots,T \label{eq:gen_limit} \\
        & P_{bus}(t) = C_g \cdot P_g(t) + C_r \cdot (P_r(t) - P_{rc}(t)) \nonumber \\
        & \hspace{1.5cm} - C_l \cdot (P_l(t) - P_{ls}(t)), \; t=1,\cdots,T \label{eq:injection} \\
        & 1^TP_{bus}(t) = 0,\; t=1,\cdots,T \label{eq:power_balance} \\
        & -P_{f,max} \leq F\cdot (P_{bus}(t) - P_{bus,s}) + P_{f,s} \leq P_{f,max}, \nonumber \\
        & \hspace{4.5cm} \; t = 1,\cdots,T \label{eq:power_flow} \\ 
        & 0 \leq P_{ls}(t) \leq P_l(t), \; t=1,\cdots,T \label{eq:slack_1}  \\
        & 0 \leq P_{rc}(t) \leq P_r(t), \; t=1,\cdots,T  \label{eq:slack_2} \\
        & \textstyle \sum_{\tau=t}^{t+T_{on}(i)-1}y_{g,i}(\tau) \leq u_{g,i}(t + T_{on}(i) -1), \nonumber \\ 
        & \hspace{0.6cm} t=1,\cdots,T - T_{on}(i) + 1, \; i=1,\cdots,N_g \label{eq:minimin_on} \\
        & \textstyle \sum_{\tau=t}^{t + T_{off}(i) -1 }z_{g,i}(\tau) \leq 1-u_{g,i}(t + T_{off}(i) -1), \nonumber \\ 
        & \hspace{0.4cm} t=1,\cdots,T - T_{off}(i) + 1, \; i=1,\cdots,N_g \label{eq:minimum_off} \\
        & \psi(u_g(t), P_r(t) - P_{tc}(t) ; \theta) \gtreqless \phi_c, \; i=1,\cdots, T \label{eq:icnn_constratint}
    \end{align}
\end{subequations}

\begin{table}[t]
\small
\centering
\caption{Parameters of the UC model.}
\begin{tabularx}{\linewidth}{|cX|}\hline
$T$ & Number of time periods \\
$N_g$ & Number of generators \\
$c_g$ & Linear generation cost coefficients (\$/MWh) \\
$c_u$ & Commitment (fixed on) cost coefficients (\$/h) \\
$c_y$ & Startup cost coefficients (\$/event) \\
$c_z$ & Shutdown cost coefficients (\$/event) \\
$c_{ls}$ & Load shedding penalty coefficients (\$/MWh) \\
$c_{rc}$ & Renewable curtailment penalty coefficient (\$/MWh) \\
$P_{min},\,P_{max}$ & Min./max.\ generator output limits (MW) \\
$R_{up},\,R_{down}$ & Ramp-up/ramp-down limits (MW/h) \\
$R_{start},\,R_{shut}$ & Startup/shutdown ramp limits (MW/h) \\
$P_l(t)$ & Bus load demand vector at time $t$ (MW) \\
$P_r(t)$ & Available renewable generation vector at time $t$ (MW) \\
$C_g$ & Generator-to-bus incidence matrix \\
$C_l$ & Load-to-bus incidence matrix \\
$C_r$ & Renewable-to-bus incidence matrix \\
$P_{bus,s}$ & Bus shift vector for transformer (MW) \\
$P_{f,s}$ & Branch flow shift vector for transformer (MW) \\
$F$ & PTDF matrix \\
$P_{f,max}$ & Line flow limits (MW) \\
$T_{on}(i)$ & Minimum on time of generator $i$ (h) \\
$T_{off}(i)$ & Minimum off time of generator $i$ (h) \\\hline
\end{tabularx}
\label{tab:param_uc_paper}
\end{table}

\begin{table}[t]
\small
\centering
\caption{Decision variables of the UC model.}
\begin{tabularx}{\linewidth}{|cX|}\hline
$P_g(t)$ & Generator dispatch vector at time $t$ (MW) \\
$u_g(t)$ & Generator on/off commitment vector at time $t$ \\
$y_g(t)$ & Generator startup indicator vector at time $t$ \\
$z_g(t)$ & Generator shutdown indicator vector at time $t$ \\
$P_{ls}(t)$ & Load shedding vector at time $t$ (MW) \\
$P_{rc}(t)$ & Renewable curtailment vector at time $t$ (MW) \\
$P_{bus}(t)$ & Net bus injection vector at time $t$ (MW) \\\hline
\end{tabularx}
\label{tab:var_uc_paper}
\end{table}

\subsubsection{ICNN as Mixed-integer Linear Constraints}\label{app:ICNN_mixed_integer}

Generic NN with linear layer and piece-wise linear activations such as ReLU can be denoted as a set of mixed-integer linear constraints \cite{xu2023availability}. Using the Big-M method, the $x=\operatorname{ReLU}(s) = \max(0,s)$ is equivalent to,
\begin{equation*}
    x\geq s, \; x\geq 0, \; x\leq u\circ \delta, \; x\leq s-l\circ (1-\delta)
\end{equation*}
where $\circ$ is the element-wise product. $x$ is post-activation output; $s$ is the pre-activation output. For linear layer, $s = W^xx_{-1} + W^zz + b$. $\delta$ is the auxiliary binary variables whose dimension equals to the output dimension of the layer. $u$ and $l$ are the bounds of the \emph{pre-activation} output $s$.

Consequently, the MILP reformulation of stability constraint optimization \eqref{eq:ed} becomes,
\begin{equation}\label{eq:MILP_stability_constraint}
    \begin{aligned}
        \min_{z,x_{1:k},\delta_{1:k-1}} \; & z^TP(y)z + q^T(y)z \\
        \text{s.t.} \; & A(y)z = b(y),\; G(y)z \leq h(y), \; x_k \lesseqgtr \phi_c \\
        & x_k = W_{k-1}^x x_{k-1} + W_{k-1}^z z + b_{k-1} \\
        & x_{i+1} \geq W_i^x x_i + W_i^z z + b_i,\; i=0,\cdots,k-2 \\
        & x_{i+1} \geq 0, \; i=0,\cdots,k-2 \\
        & x_{i+1} \leq u_{i+1} \circ \delta_{i+1}, \; i=0,\cdots,k-2 \\
        & x_{i+1} \leq W_{i}^x x_i + W_i^z z + b_i - l_{i+1} \circ (1-\delta_{i+1}),  \\
        & \hspace{3.5cm} i=0,\cdots,k-2
    \end{aligned}
\end{equation}

To fairly compare with the mixed-integer reformulation of ICNN as stability constraint, the interval bound propagation (IBP) is implemented, which computes valid lower and upper bounds ($\tilde{l}$ and $\tilde{u}$) of the post-activation of each layer. It will be used to determine the bounds of pre-activation of each layer, which is essential in big-M method \cite{gowal2018effectiveness}. Although tighter approximations exist, the IBP in general provides the big-M sufficient for solving an optimization problem \cite{xu2023availability}. 
\begin{equation*}
    \begin{aligned}
        {l}_{i+1} & = \max(W_i^x,0) \cdot \tilde{l}_i + \min(W_i^x,0) \cdot \tilde{u}_i + \max(W_i^z,0) \cdot \underline{z} \\& \hspace{1.5cm} + \min(W_i^z,0) \cdot \bar{z} + b_i, \; i=0,\cdots,k-2 \\
        {u}_{i+1} & = \max(W_i^x,0) \cdot \tilde{u}_i + \min(W_i^x,0) \cdot \tilde{l}_i + \max(W_i^z,0) \cdot \bar{z} \\& \hspace{1.5cm} + \min(W_i^z,0) \cdot \underline{z} + b_i, \; i=0,\cdots,k-2 \\
        \tilde{l}_{i+1} & = \operatorname{ReLU}({l}_{i+1}), \; \tilde{u}_{i+1} = \operatorname{ReLU}({u}_{i+1}), \; i=0,\cdots,k-2
    \end{aligned}
\end{equation*}
where $\underline{z} \leq z \leq \bar{z} $ are the lower and upper bounds of the decision variables. In the simulation setting on \eqref{eq:ed_over_time}, $\underline{z}$ is a vector of zero, representing the both the minimum of $u_g$ and $P_{rn}$ are zero. Similarly, $\bar{z}$ is a vector of one and the available renewable generation at the moment $P_r(t),\; t=1,\cdots,T$.

\subsubsection{Outer Approximation Algorithm}\label{app:outer_algorithm}

\begin{figure}[t]
    \centering
    \includegraphics[width=0.8\linewidth,trim={5 5 5 5 cm},clip]{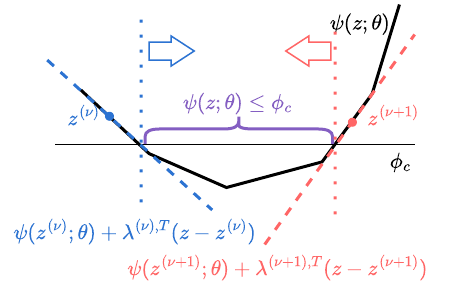}
    \caption{Illustration on the OA cuts for single time step case.}
    \label{fig:oa_cut}
\end{figure}

For outer approximation, \eqref{eq:ed_over_time} is set with ``$\leq$''. In this case, the label ``0'' should be marked as the stable samples when constructing the dataset. To ensure stability on unseen scenarios, $\phi_c$ is set to be slightly smaller than 0. The following standard outer approximation is modified from \cite{conejo2006decomposition}.

\textbf{Step 0: Initialization.} Initialize the iteration counter $\nu=0$ and solve \eqref{eq:ed_over_time} with no ICNN surrogate constraint as $z^{(0)}$.

\textbf{Step 1: Determining constraint violations.} Identify the most violated ICNN surrogate constraint, with time index $t^\star_{\nu}$. 
\begin{equation}
    \psi(z_{t^\star_{\nu}}^{(\nu)};\theta) = \max_{t=1:T} \psi(z_{t}^{(\nu)};\theta)
\end{equation}
If $\psi(z_{t^\star_\nu}^{(\nu)};\theta) \leq \phi_c$, stop and return $z^{(\nu)}$; otherwise, continue to Step 2.

\textbf{Step 2: Linearization.} Obtain the subgradient of the \emph{most} violated ICNN surrogate with respect to $z^{(\nu)}_{t^\star_\nu}$, i.e., $\lambda^{(\nu)} \in \partial \psi(z^{(\nu)}_{t^\star};\theta)$. This step can be achieved by automatic differentiation package such as \texttt{PyTorch} \cite{paszke2019pytorch} or by solving \eqref{eq:lp_equivalence_subproblem}. A new cut is constructed as
\begin{equation}
    \ell^{(\nu)}(z_{t^\star_{\nu}}) = \psi(z^{(\nu)}_{t^\star_\nu};\theta) + \lambda^{(\nu),T} (z_{t^\star_\nu} - z^{(\nu)}_{t^\star_\nu})
\end{equation}

It is also possible to use the LP representation of the ICNN in \eqref{eq:lp_equivalence} to compute subgradients. In particular, one may solve the following subproblem with $z$ fixed at the current iterate:
\begin{equation}\label{eq:lp_equivalence_subproblem}
    \begin{aligned}
        & \psi(z^{(\nu)};\theta) := \textstyle \min_{x_{1:k}} x_k \\
        \text{s.t.} \; & x_{i+1} \geq W_i^x x_i + W_i^z z^{(\nu)} + b_i,\; i=0,\cdots,k-1 \\
        & x_i \geq 0,\; i=1,\cdots,k-1.
    \end{aligned}
\end{equation}
Let $\mu^{(\nu)}$ denote the optimal dual multipliers associated with the constraints $x_{i+1} \geq W_i^x x_i + W_i^z z^{(\nu)} + b_i$. Then a subgradient can be obtained from the dual variables (equivalently, from sensitivity of the optimal value with respect to $z$). 

\textbf{Step 3: Solution on the linearized problem.} Solve the following problem,
\begin{equation}\label{eq:ed_over_time_outer}
    \begin{aligned}
         \min_{z_{1:T}} \; & z^TP(y)z + q^T(y)z \\
        \text{s.t.} \; & A(y)z = b(y),\; G(y)z \leq h(y) \\
        & \ell^{(s)}(z_{t^\star_\nu}) \leq \phi_c,\; s = 1,\cdots,\nu
    \end{aligned}
\end{equation}
The solution is obtained as $z^{(\nu+1)}$. Set $\nu\leftarrow\nu+1$ and return \textbf{Step 1}.

The OA algorithm is illustrated in Fig.~\ref{fig:oa_cut} for single time-step case. After the two successive cuts, the feasible region $\psi(z;\theta)\leq\phi_c$ is recovered. Moreover, it also visualized that for piece-wise convex function $\psi(z;\theta)$, the global convergence can be attained, as proved in Appendix~\ref{app:proof_oa}.  Meanwhile it is possible to take \emph{all} violated constraints as new cuts, which is followed in the experiment.

\subsubsection{Inner Approximation}\label{app:inner_algorithm}

\begin{figure}[t]
    \centering
    \subfloat[Iteration $\nu$]{%
        \includegraphics[width=0.7\linewidth,trim={0 5 5 5 cm},clip]{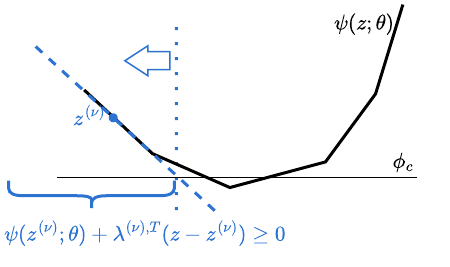}
        \label{fig:ia_cut_step_one}
    }
    \hfill
    \subfloat[Iteration $\nu+1$]{%
        \includegraphics[width=0.7\linewidth,trim={0 5 5 5 cm},clip]{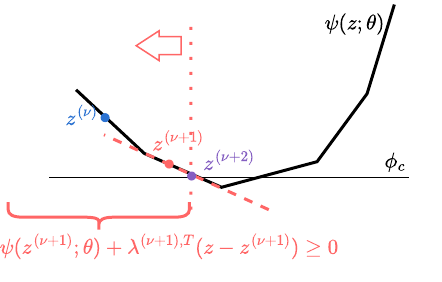}
        \label{fig:ia_cut_step_two}
    }
    \caption{Illustration on the IA cuts for single time step case.}
    \label{fig:ia_cut}
\end{figure}

For the inner approximation setting, the stable label is set as ``1'' so that the trained ICNN is constrained as $\psi(z,\theta) \geq \phi_c$ to ensure stability. Similarly as the outer approximation case, $\phi_c$ is set to be slightly larger than 0. 

Unlike outer approximation, the initial point $z^{(0)}$ should be feasible for \eqref{eq:ed_over_time}. A two-phase procedure is implemented where the Phase I is to ensure the feasibility and Phase II is for optimality. 

\underline{\emph{Phase-I Problem:}}

\textbf{Step 0: Initialization.} Initialize the iteration counter $\nu=0$ and solve \eqref{eq:ed_over_time} with no ICNN surrogate constraint as $z^{(0)}$.

\textbf{Step 1: Linearization.} Obtain the subgradient $\lambda^{(\nu)}_t \in \partial \psi(z^{(\nu)}_t;\theta),\; t=1,\cdots,T$ and construct the following cuts,
\begin{equation}
    \ell^{(\nu)}(z_{t}) = \psi(z^{(\nu)}_{t};\theta) + \lambda^{(\nu),T}_t (z_{t} - z^{(\nu)}_{t}), \; t=1,\cdots,T
\end{equation}

\textbf{Step 2: Solution on the linearized problem.} Solve the following optimization problem.
\begin{equation}
    \begin{aligned}
        \min_{s_{1:T},z_{1:T}} \; & \textstyle \sum_{t=1}^T s_t + \frac{\rho}{2} (s_t - s^{(\nu-1)}_t )^2 \\
        \text{s.t.} \; & A(y)z = b(y),\; G(y)z \leq h(y) \\
        & \ell^{(\nu)}(z_t) + s_t \geq \phi_c, \; t=1,\cdots,T \\
        & s_t \geq 0,\; t=1,\cdots,T
    \end{aligned}
\end{equation}
where $s_t,\;t=1,\cdots,T$ are the slack variables for constraint violation. The solution is denoted $z^{(\nu+1)}$. If $\sum_{t=1}^T s_t = 0$, stop and return $z^{(\nu+1)}$ as $z^{(0)}$ in Phase II; otherwise, set $\nu\leftarrow\nu+1$ and return to \textbf{Phase I-Step 1}.

\underline{\emph{Phase-II Problem:}}

\textbf{Step 0: Initialization.} Initialize $\nu=0$ and obtain the feasible $z^{(0)}$ from Phase 1.

\textbf{Step 1: Linearization.} Obtain the set of cuts as in Phase I-Step 1.

\textbf{Step 2: Solution on the linearized problem.} Solve the following optimization problem.
\begin{equation}
    \begin{aligned}
        \min_{z_{1:T}} \; & z^TP(y)z + q^T(y)z + \textstyle \frac{\rho}{2} \sum_{t=1}^T\|z_t - z_t^{(\nu-1)}\|^2 \\
        \text{s.t.} \; & A(y)z = b(y),\; G(y)z \leq h(y) \\
        & \ell^{(\nu)}(z_t) \geq \phi_c, \; t=1,\cdots,T \\
    \end{aligned}
\end{equation}
The solution is $z^{(\nu+1)}$. If the optimal value between $\nu$ and $\nu+1$ is smaller than a tolerance, stop and return $z^{(\nu+1)}$; otherwise, set $\nu\leftarrow \nu+1$ and return to \textbf{Step 1}.

Two successive cuts are visualized in Fig.~\ref{fig:ia_cut}. Unlike the OA algorithm, new IA cuts \emph{replace} the cuts from previous iterations. Moreover, only the left side of the convex function $\psi(z;\theta)$, where initial point $z^{(\nu)}$ locates, is explored. Therefore, no global optimum is guaranteed.

\subsubsection{Discussion on the Outer and Inner Approximations}\label{app:outer_vs_inner}

As shown by Table~\ref{tab:compare_outer_inner}, the dual view between the two approximation methods become clear. For OA, each iteration step adds ICNN surrogate constraint as new linearized cut while \emph{keeps} the previous cuts. The cumulative cuts become the lower bound of the convex ICNN surrogate. In contrast, the IA adds linearized cuts for \emph{all} ICNN surrogate constraints over $t=1,\cdots,T$ while the previous cuts are \emph{removed}.

\begin{table}[t]
\centering
\caption{Comparison between outer and inner approximations}
\label{tab:compare_outer_inner}
\footnotesize
\setlength{\tabcolsep}{4pt}
\renewcommand{\arraystretch}{1.1}
\begin{tabularx}{\columnwidth}{lXX}
\toprule
 & \textbf{Outer Approx.} & \textbf{Inner Approx.} \\
\midrule
\textbf{Convergence}
& Guaranteed
& Not guaranteed \\

\textbf{1st iter. Feasibility}
& Not guaranteed
& Guaranteed \\

\textbf{Cut behavior}
& Add cut(s) cumulatively
& Replace cuts for all surrogate constraints \\
\bottomrule
\end{tabularx}
\end{table}


\subsubsection{Optimization Structures}\label{app:optimization_structure}

The exact number of extra constraints and variables compared to the original UC are summarized in Table~\ref{tab:extra_component}, which can be used to verify the result in Table~\ref{tab:14_bus} and \ref{tab:118_bus}.

\begin{table}[t]
    \centering
    \setlength{\tabcolsep}{3pt}
    \renewcommand{\arraystretch}{0.9}
    \caption{Optimization structure of ICNN-based small-signal stability constrained optimization. $T$ (=24) is time horizon; $N_H$ is the total hidden dimension in ICNN; $N_I$ (=5 for 14-Bus system and =54 for 118-Bus system) is the ICNN input size. $R$ is the number of iteration in outer approximation. $N_\nu$ is the number of violated constraints in the $\nu$-th iteration.}
    \begin{tabular}{ccc}
        \toprule
         &  \textbf{No. Extra Cons.} & \textbf{No. Extra Vars., (Binary)}  \\
         \midrule
       D-MILP-$\geq$/$\leq$  & $T\times(N_I + 4N_H+2)$ & $T \times (N_I + 2N_H+1)$, ($T\times N_H$)  \\
       D-LP-$\leq$ & $T\times (N_I + 2N_H + 2)$ & $T\times(N_I + N_H + 1)$, ($0$) \\
       I-OA-$\leq$ & $\sum_{\nu=1}^{R} N_\nu$ & $0$, ($0$) \\
       I-IA-$\geq$ & $T$ & $0$, ($0$) \\
       \bottomrule
    \end{tabular}
    \label{tab:extra_component}
\end{table}

\subsection{Simulation Settings}\label{app:simulation_setting}

We use workstation with two Intel(R) Xeon(R) Gold 6230R CPU @ 2.10GHz and up to 104 cores with Ubuntu 20.04.6 LTS. This is a relatively powerful computational unit for medium to large scale power system optimizations. The ICNN is trained with NVIDIA GeForce RTX 3090. The hyperparameters for ICNN structure and training are reported in Table~\ref{tab:icnn_config}. All the optimization problems are solved by \texttt{Gurobi} with \texttt{MIPGap}=1e-3, \texttt{MIPGapAbs} = 1e-3, \texttt{MIPFocus}=1, \texttt{Heuristics}=0.5, and \texttt{TimeLimit} = 3600.

For 14-Bus system, two SGs at bus 3 and 6 are replaced by solar power plant with the same rated capacities. 21 out of 54 SGs are replaced in 118-Bus system. Other configurations can be found on GitHub repository. Active sampling strategy is implemented to augment the dataset constructed by running the basic UC without stability constraints on realistic load and solar generation profiles. The raw dataset is obtained from Texas power system \cite{lu2025synthetic} and preprocessed to match the scalability of 14 and 118 Bus systems by our open-source package \texttt{GridForge} \cite{xu2025lapso}. The critical gSCR value is set as 2.5 and we also slightly increase the threshold during data generation to increase the true positive rate. 

Active sampling strategy in \cite{xu2025lapso} is adopted to build the dataset. Given the one-year load and solar generation profile, operational points $\{(u_g(t),P_{rn}(t))\}_{t=1}^{8760}$ are solved via basic operation without stability constraints. For each operational point with label $y(t)$, the closest operational point with \emph{opposite} stability label is iteratively obtained through gradient-based and heuristic-based sampling method (See \cite{xu2025lapso} for details). The active sampling strategy augments the dataset to be close to the decision boundary, which are informative for stability constraint learning. The size of the resultant dataset is 78{,}840, which is denoted as $\mathcal{D} = \{(u_g(t),P_{rn}(t)),y(t)\}_{t=1}^{78840}$ where $y(t)\in\{0,1\}$. The choice of positive label depends on either the ``$\leq$'' and ``$\geq$'' settings.

\begin{table}[h]
    \centering
    \caption{ICNN Configurations}
    \begin{tabular}{c|c|c|c}
    \toprule
       \textbf{Hidden size}  & (50, 20) & \textbf{Learning rate} & 1e-3 \\
       \textbf{Optimizer} & ADAM & \textbf{Hidden activation} & ReLU \\
       \textbf{Batch size}  & 512 & \textbf{Epochs} & 100 \\
       \textbf{Learning rate decay} & 0.8 at epoch 50 & \textbf{Early stop} & 10 \\
       \bottomrule
    \end{tabular}
    \label{tab:icnn_config}
\end{table}

\renewcommand{\theequation}{S-\arabic{equation}}
\setcounter{equation}{0}

\renewcommand{\thetable}{S-\arabic{table}}
\setcounter{table}{0}

\renewcommand{\thefigure}{S-\arabic{figure}}
\setcounter{figure}{0}

\bibliographystyle{IEEEtran}
\bibliography{IEEEabrv,Reference.bib}


 




\vfill

\end{document}